\documentclass[12pt]{article}
\usepackage[dvips]{graphicx}
\usepackage{amsmath,amssymb}
\usepackage{bm}

\title{On NP complete problems I}
\author{}

\begin{document}
{
\rightline{September 2008}
\rightline{~~~~~~~~~}
\vskip 1cm
\centerline{\Large \bf On NP complete problems I}
\vskip 1cm
\centerline{Minoru Fujimoto${}^{a)}$ and Kunihiko Uehara${}^{b)}$}
\vskip 1cm
\centerline{${}^{a)}${\it Seika Science Research Laboratory,
Seika-cho, Kyoto 619-0237, Japan}}
\centerline{${}^{b)}${\it Department of Physics, Tezukayama University,
Nara 631-8501, Japan}}
\vskip 5mm
\centerline{\bf Abstract}
\vskip 0.2in
  We study the quadratic residue problem known as an NP complete problem 
by way of the prime number and show that a nondeterministic polynomial process 
does not belong to the class P because of a random distribution of solutions 
for the quadratic residue problem.


\vskip 12mm
}
\setcounter{equation}{0}
\addtocounter{section}{0}
\section{Introduction}
\hspace{\parindent}

  It is known that there is ``a quadratic residue problem" or ``a quadratic 
congruence problem" in the NP complete problems\cite{Manders}, 
which is a problem that ``Is there a positive integer $x$ whose value is 
less than $c$ and is such that 
\begin{equation}
  x^2\equiv a\ (\text{mod } b),
\label{e101}
\end{equation}
{\it i.e.} the remainder when $x^2$ is divided by $b$ is equal to $a$?" 

  We will show that a special case of this quadratic residue problem 
does not belong to the class P, namely, 
it will be proved that there is no skillful method to determine a solution 
in polynomial time for a nondeterministic portion which will be 
essentially a super-polynomial process in the NP complete problem. 

Generally speaking, by the Fermat's little theorem 
\begin{equation}
  n^p\equiv n\ (\text{mod } p),
\label{e102}
\end{equation}
a polynomial, whose coefficient is made of a natural number 
$n$ to the some power multiplied by an integer, 
turns into a periodic function in modulo $p$ for a variable $x=kp$. 
More specifically, a function $f(x)$, 
whose coefficients are rational and constant term is integer, 
becomes a periodic function in modulo $p$ for any value of the argument, 
where $p$ is any prime number, 
a form of an argument is $kp\times(\text{lcm of $d$'s})+r$ with that
$d$'s are denominators of the rational coefficients of $f(x)$ and 
$r$ is an integer between $0$ and $p-1$.

Therefore, it should be proved that the time required for 
the decision of a nondeterministic process must be expressed 
in a super-polynomial time for a modulo operation.

\section{A Proof}
\hspace{\parindent}

  Now we confine ourselves to a special case that $b$ is a prime number of 
$p=4m+1$ type in Eq.(\ref{e101}).
  We take $b$ as $p=4m+1$, 
then $p$ is uniquely expressed by 
\begin{equation}
  p=s^2+t^2,
\label{e201}
\end{equation}
where $s$ and $t$ are natural numbers with $s>t$, 
which is shown by Fermat. 
The another expression
\begin{equation}
  s^2\equiv p-t^2\ (\text{mod } p)
\label{e202}
\end{equation}
is easily got by taking modulo $p$ in Eq.(\ref{e201}), 
which is the same form of the quadratic residue problem 
for $x=s$, $a=p-t^2$ and $b=p$ in Eq.(\ref{e101}). 

  On the other hand, by the Wilson's theorem\cite{Mordell,Chowla}
\begin{equation}
  (p-1)!\equiv -1\ (\text{mod } p)
\label{e203}
\end{equation}
is satisfied for any prime which can be read as 
\begin{equation}
  \left\{\left(\frac{p-1}{2}\right)!\right\}^2\equiv -1\ (\text{mod } p)
\label{e204}
\end{equation}
for primes of $4m+1$ type. 
So the general solution for Eq.(\ref{e101}) is given by 
$\displaystyle{x=\left(\frac{p-1}{2}\right)!\ (\text{mod }p)}$ 
with $a=p-1$ and $b=p$, 
and this is the case of $t=1$ in Eq.(\ref{e201}) where we easily check 
whether $s$ or $x$ is less than $c$ if $s$ is an integer. 
But generally we have to search a solution $s$ satisfying Eq.(\ref{e202}) 
for a given $a$. 

  What we have to show is whether an amount of the computational complexity 
for fixing $s$ is over a deterministic polynomial time. 
  As is stated above, there is only one solution $s$ and $t$ 
for a given $p$, then it can be said that $s$ and $t$ distribute at random 
as follows.

 The equation (\ref{e201}) can be regarded as a product of Gaussian primes, 
so we refer the extended prime number theorem for the Gaussian primes as
\begin{equation}
  \pi(x,\mathbb{Z}[i])=\sharp\left\{\alpha\in\mathbb{Z}[i]
      \text{prime elements}
  \bigg|
    \begin{matrix}
       |\alpha|^2\le x\\
       0\le\arg(\alpha)<\pi/2
    \end{matrix}
  \right\}\sim \frac{x}{\log x}
\label{e207}
\end{equation}
and it can be said that there distribute uniformly for arguments with 
$0\le \theta_1<\theta_2<\frac{\pi}{2}$,
\begin{equation}
  \pi(x,\mathbb{Z}[i];\theta_1,\theta_2)
  =\sharp\left\{\alpha\in\mathbb{Z}[i]
      \text{prime elements}
  \bigg|
    \begin{matrix}
       |\alpha|^2\le x\\
       \theta_1\le\arg(\alpha)\le\theta_2
    \end{matrix}
  \right\}\sim \frac{2}{\pi}(\theta_1-\theta_2)\frac{x}{\log x}.
\label{e208}
\end{equation}
The theorem makes us to be confident that $s$ and $t$ distribute at random. 
This means that $s$ and $t$ are elements of non-polynomial Gaussian integers. 
This coincides the fact that there is a contradiction that $p$ can be 
expressed by a finite polynomial, if $s,t$ are expressed by polynomials.
Even in the case of $t=1$, solutions $s$ for Eq.(\ref{e202}) 
distribute at random because of the modulo $p$. 

  When the solutions for Eq.(\ref{e202}) distribute at random, we have to 
check every number less then $\displaystyle{\sqrt{\frac{p-1}{2}}}$ 
for a candidate of $x$
satisfying the equation whether there exits 
$x$ less than $c$. 
This makes us to increase an amount of the computational complexity 
and we need 
over a deterministic polynomial time. 

  Because any other NP complete problem can be reduce to the quadratic residue 
problem in polynomial time, 
it is possible to prove that all NP complete problems can not be solved 
in a polynomial calculated amount.
Although a primality test in a polynomial amount exists such as 
the AKS method, a solution for a factorial residue of primes 
needs calculations over a polynomial amount. 

  Above all, there exists no deterministic method for the NP complete problem 
in polynomial time because of a random distribution of the solutions 
for the quadratic residue problem.

\section{Conclusions}
\hspace{\parindent}

@As is stated in previous section, the general solution $x$ 
in Eq.(\ref{e101}) is expressed by a factorial, 
whose generalization is the gamma function. 
If a value of the gamma function is got in polynomial time, we easily judge 
an integer whether it is a prime number or not using the Wilson's theorem. 
But the reason why we had to wait for the AKS method is that this is 
not the case. 

  The gamma function itself is non-polynomial, which is understood by a fact 
that there is no differential equation for the gamma function, 
and the gamma function modulo $p$ is a non-periodic function, 
so it is impossible to display in polynomial time. 

  As is shown in Appendix, an essential point of the quadratic residue 
problem is that there is no skillful method to count up lattice points on 
a circle with a radius $r$ when $r^2$ is a prime number except $2$. 
  That is why the prime above is a type of $4m+1$ and $s,t$ satisfying 
$s^2+t^2=p$ distribute at random. 
So we have to check every number less then 
$\displaystyle{\sqrt{\frac{p-1}{2}}}$ for $s$, 
this makes us to increase an amount of the computational complexity.

\vskip 15mm
\renewcommand{\theequation}{\Alph{section}\arabic{equation}}
\setcounter{section}{1}
\setcounter{equation}{0}
\section*{Appendix}
\hspace{\parindent}
  Here we give another view of that $s$ and $t$ satisfying 
$s^2+t^2=p$ distribute at random. 

The number of lattice points $x,y$ on or inside a circle with a radius $r$ 
is $[\pi r^2]$, where $[\ ]$ denotes a Gauss symbol, 
whereas the number of lattice points $N_0$ for $x\ge y>0$ is around
$[\pi r^2/8]-r$, 
so we can write down the number as $N_0=\pi r^2/8-cr+o(r)$

When we put $r^2=R$, then
$$
  N_0=\frac{\pi R}{8}-c\sqrt{R}+o(\sqrt{R}).
$$
Primes corresponding to lattice points up to $R$ are $2$ and 
$4m+1$ type primes, their number $N$is given by the prime number theorem as
$$
 N\sim \frac{\pi(R)}{2}+1\sim \frac{R}{2\log R},
$$
$$
 \frac{N}{N_0}\sim \frac{4}{\pi\log R},
$$

 @Lattice points corresponding to primes distribute at random 
because lattice points distribute uniform but primes of $4m+1$ type distribute 
at random, 
so the lattice points of solution $s,t$ distribute at random 
over the $\theta$ direction in the polar coordinates.

\vskip 10mm

\end{document}